\documentclass{amsart}
\usepackage{amssymb}

\newcommand{\ZFCa}{{\operatorname{\mathsf {ZFC}}}}
\newcommand{\CH}{\operatorname{\mathsf {CH}}}

\newcommand{\reals}{{\mathbb R}}

\newcommand{\rest}{{\mathord{\restriction}}}
\newcommand{\add}{\operatorname{\mathsf   {add}}}
\newcommand{\cov}{\operatorname{\mathsf  {cov}}}
\newcommand{\unif}{\operatorname{\mathsf  {non}}}

\newcommand{\cf}{{\operatorname{\mathsf   {cf}}}}

\newcommand{\dom}{{\operatorname{\mathsf {dom}}}}

\newcommand{\supp}{{\operatorname{\mathsf    {supp}}}}

\newcommand{\QED}{\hspace{0.1in} \square \vspace{0.1in}}
\newcommand{\forces}{\Vdash}

\newcommand{\N}{{\mathcal N}}
\newcommand{\M}{{\mathcal M}}
\newcommand{\V}{{\mathbf V}}

\newcommand{\thinks}{\models}

\newcommand{\Proof}{{\sc Proof} \hspace{0.2in}}

\newcommand{\lft}[2]{\mathopen\ifcase#1{}\oo\or
                        \big#2\or\Big#2\else\oo\fi} 
\newcommand{\rgt}[2]{\mathclose\ifcase#1{}\oo\or
                        \big#2\or\Big#2\else\oo\fi}

\newcommand{\SM}{{\mathcal{SM}}}
\newcommand{\SN}{{\mathcal  {SN}}}

\theoremstyle{plain}
\newtheorem{theorem}{Theorem}[section]
\theoremstyle{plain}
\newtheorem{lemma}[theorem]{Lemma}

\newtheorem{definition}[theorem]{Definition}

\begin{document}

\title{Strongly meager and strong measure zero sets}
\author{Tomek Bartoszy\'{n}ski}
\address{Department of Mathematics and Computer Science\\
Boise State University\\
Boise, Idaho 83725 U.S.A.}
\thanks{First author supported by Alexander von Humboldt Foundation and
NSF grant DMS 95-05375}
\email{tomek@math.idbsu.edu, http://math.idbsu.edu/\char 126 tomek}
\author{Saharon Shelah}
\thanks{Second author partially supported by Basic Research Fund,
Israel Academy of Sciences, publication 658}
\address{Department of Mathematics\\
Hebrew University\\
Jerusalem, Israel}
\email{shelah@sunrise.huji.ac.il, http://math.rutgers.edu/\char 126 shelah/}

\begin{abstract}
In this paper we present two consistency results concerning
the existence  of large strong measure
zero and strongly meager sets.
\end{abstract}
\maketitle
\section{Introduction}

Let $\M$ denote the collection of all 
meager subsets of $2^\omega $ and let $\N$ be the collection of all subsets of $2^\omega$ that have measure zero with respect to 
the standard product measure on $2^\omega $.

\begin{definition}
  Suppose that $X \subseteq 2^\omega $ and let $+$ denote the
componentwise addition modulo $2$.  We say that $X$ is strongly meager
if for every $H \in \N$, $X+H =\{x+h:x \in X, h \in H\}\neq 2^\omega
$. 

We say that $X$  is a strong measure zero set if for every $F \in
  \M$, $X+F \neq 2^\omega $. Let $\SM$ denote the collection of strongly meager sets and let $\SN$ denote the collection of
  strong measure zero sets.

For a family  of sets ${\mathcal J} \subseteq P(\reals)$ let

$\cov({\mathcal J}) = \min\left\{ |{\mathcal A}| : {\mathcal A} \subseteq {\mathcal J}
\hbox{ and } \bigcup {\mathcal A} = 2^\omega\right\}$.

$\unif({\mathcal J}) = \min\left\{ |X| : X \not\in {\mathcal J} \right\}$.
\end{definition}
Strong measure zero sets are usually defined as those subsets $X$ of
$2^\omega$   such that for every sequence of positive reals
$\{\varepsilon_n: n \in \omega\}$ there exists a sequence of basic
open sets $\{I_n:n \in\omega\}$ with diameter of $I_n$ smaller than
$\varepsilon_n$ and $X \subseteq \bigcup_n I_n$.
The Galvin-Mycielski-Solovay theorem (\cite{GaMySo}) guarantees
that both definitions are yield the same families of sets.

Recall the following well--known facts.
Any of the following  sentences is consistent with $\ZFCa$,
\begin{enumerate}
\item $\SN=[2^\omega]^{\leq
    \boldsymbol\aleph_0}$, (Laver \cite{Lav76Con})
\item $\SN=[2^\omega]^{\leq
    \boldsymbol\aleph_1}$, (Corazza \cite{Cor89Gen}, Goldstern-Judah-Shelah \cite{GJS92})
\item  $\SM=[2^\omega]^{\leq
    \boldsymbol\aleph_0}$. (Carlson, \cite{CarStr})
\item   $\unif(\SN)= {\mathfrak d}=2^{\boldsymbol\aleph_0}> \boldsymbol\aleph_1 $,
$\cov(\M)=\boldsymbol\aleph_1 $ and there exists a strong measure zero
set of size $ 2^{\boldsymbol\aleph_0} $. (Goldstern-Judah-Shelah \cite{GJS92})
\end{enumerate}
The proofs of the above results as well as all other results quoted in
this paper can be also found in 
\cite{BJbook}. 

In this paper we will show that the following statements are consistent
with $\ZFCa$:
\begin{itemize}
\item for any regular $ \kappa >
  \boldsymbol\aleph_0 $, $\SM=[2^\omega]^{< \kappa }$,
\item $\SM$ is an ideal and $\add(\SM) \geq \add(\M)$, 
\item  $\unif(\SN)=2^{\boldsymbol\aleph_0}>
  \boldsymbol\aleph_1  $, ${\mathfrak d}=\boldsymbol\aleph_1$ and
  there is a strong measure zero set of size $ 2^{\boldsymbol\aleph_0}
  $.
\end{itemize}

\section{$\SM$ may have large additivity}
In this section we will show that $\SM$ can be an ideal with large
additivity. 
Let 
$${\mathfrak m}=\min\{\gamma: {\mathbf {MA}}_\gamma \text{ fails}\}.$$
We will show that $\SM=[2^\omega]^{<{\mathfrak m}}$ is consistent with
$\ZFCa$, provided ${\mathfrak m}$ is regular. In particular, the model
that we construct will satisfy
$\add(\SM)=\add(\M)$.

Note that if $\SM=[2^\omega]^{<{\mathfrak m}}$ then $
2^{\boldsymbol\aleph_0} > {\mathfrak m}$, since Martin's Axiom implies
the existence of a strongly meager set of size $2^{\boldsymbol\aleph_0}$.
Our construction is a generalization of the construction from \cite{CarStr}.

To witness that a set is not strongly meager we need a measure zero
set. The following theorem is crucial.

\begin{theorem}[Lorentz]\label{lorenz}
  There exists a function $K \in \omega^\reals $ such that for every 
$ \varepsilon > 0 $,
if $A \in [2^\omega]^{\geq K(\varepsilon)}$ then for all except
finitely many $k \in \omega$ there 
exists  $C \subseteq 2^k$ such
that
\begin{enumerate}
\item $|C|\cdot 2^{-k} \leq \varepsilon $,
\item $(A\rest k) +C =2^k$.
\end{enumerate}
\end{theorem}
\Proof
Proof of this lemma can be found in \cite{loradd54}
or  \cite{BJbook}.~$\QED$

\begin{definition}\label{special}
For each $n\in\omega$ let $\{C^n_m: n,m \in \omega\}$ be an enumeration of all clopen sets
in $2^\omega $ of measure
$\leq 2^{-n}$.
For a real $r \in \omega^\omega $ and $ n \in \omega $ define an open
set
$$H_n^r = \bigcup_{m >n} C^m_{r(m)}.$$
It is clear that $H_n^r$ is an open set of measure not exceeding
$2^{-n}$. In particular,
$H^r = \bigcap_{n \in \omega} H_n^r$ is a
Borel measure zero set of type $G_\delta $.

\end{definition}

\begin{theorem}
  Let $ \kappa > \boldsymbol\aleph_0 $ be a regular cardinal. It is
  consistent with $\ZFCa$ that $\mathbf{MA}_{<\kappa} +
  \SM=[2^\omega]^{<\kappa}$ holds.  In particular, it is consistent
  that $\SM $ is an ideal and
  $\add(\SM)=\add(\M)>\boldsymbol\aleph_1$.
\end{theorem}
\Proof
Fix $ \kappa $ such that $ \cf(\kappa)= \kappa >
\boldsymbol\aleph_0 $.  Let $ \lambda > \kappa $ be a regular
cardinal such that $ \lambda^{< \lambda } = \lambda$. Start
with a model $\V \thinks \ZFCa + 2^{\boldsymbol\aleph_0} =
\lambda$.

Suppose that $ {\mathcal P} $ is a forcing notion of size $< \kappa $.
We can assume that there is $ \gamma < \kappa $ such that $
{\mathcal P} = \gamma $ and $ \leq, \perp \subseteq \gamma\times
\gamma $.

Let $\{ {\mathcal P}_\alpha, \dot{{\mathcal Q}}_\alpha: \alpha<
\lambda \}$ be a finite support iteration 
such that 
for each $ \alpha < \lambda $,
\begin{enumerate}
\item $\forces_\alpha \dot{{\mathcal Q}}_\alpha \simeq {\mathbf C}$,
  if $ \alpha$ is limit,
\item there is $\gamma=\gamma_\alpha$  such that 
$\forces_{\alpha} \dot{{\mathcal Q}}_\alpha \simeq (\gamma, \leq
  , \perp)$ is a ccc forcing
  notion.
\end{enumerate}
By passing to a dense subset  we can assume that if $p \in  {\mathcal
  P}_\lambda $ 
then $p: \dom(p) \longrightarrow \kappa $, where $\dom(p)$ is a finite
subset of $ \lambda $.

By
bookkeeping we can guarantee that $\V^{{\mathcal P}_\lambda} \thinks
{\mathbf 
  {MA}}_{<\kappa}$. 
In particular, $\V^{{\mathcal P}_\lambda} \thinks [2^\omega]^{<\kappa}
\subseteq \SM$.

It remains to show that no set of size $ \kappa $ is strongly meager.

Suppose that $X \subseteq \V^{{\mathcal P}_\lambda} \cap 2^\omega $ is
a set of size $ \kappa $. Find limit ordinal  $ \alpha < \lambda $ such that $X
\subseteq 2^\omega \cap \V^{{\mathcal P}_\alpha}$. As usual we can
assume that $\alpha=0$.
Let $c$ be the Cohen real added at the step $ \alpha=0 $. 
We will show that $\V^{{\mathcal P}_\lambda} \thinks X+H^c=2^\omega $,
which will end the proof.

Suppose that the above
assertion is false. 
Let $p \in {\mathcal P}_\lambda $ and let $ \dot{z}$ be a $ {\mathcal
  P}_\lambda $-name for a real such that 
$$p \forces_\lambda \dot{z} \not \in X+H^{\dot{c}}.$$
Let $X=\{x_\xi: \xi<\kappa\}$ and for each $\xi$ find $p_\xi \geq p$
and $n_\xi \in \omega $ such that 
$$p_\xi \forces_\lambda \dot{z} \not \in x_\xi +
H^{\dot{c}}_{n_\xi}.$$

Let $Y \subseteq \kappa $ be a set of size $ \kappa $ such that
\begin{enumerate}
\item $n_\xi = \widetilde{n}$ for $\xi \in Y$,
\item $\{\dom(p_\xi): \xi \in Y\}$ form a $\Delta$-system with root
  $\widetilde{\Delta}$,
\item $p_\xi \rest \widetilde{\Delta} = \widetilde{p}$, for $ \xi \in Y$,
\item $p_\xi(0)=\widetilde{s}$, with $|\widetilde{s}|=\ell>\widetilde{n}$,
 for $\xi \in Y$.
\end{enumerate}

Fix a subset $X'=\{x_{\xi_j} : j < K(2^{-\ell})\} \subseteq Y$ and let
$\widetilde{m} \in \omega$ be such that $C^{\ell}_{\tilde{m}}+X'=2^\omega $.

Define condition $p^\star$ as
$$p^\star(\beta)=\left\{
  \begin{array}{ll}
p_{\xi_j} & \text{if } \alpha \neq \beta \ \&\  \beta \in \dom(p_{\xi_j}), \ j <K(2^{-\ell})\\
\widetilde{s}^\frown \widetilde{m} & \text{if } \alpha = \beta 
  \end{array}\right. \text{ for }\beta<\lambda .$$
On one hand $p^\star \forces_\lambda C^{\ell}_{\tilde{m}} \subseteq
H^{\dot{c}}_{\tilde{n}}$, so $p^\star \forces_\lambda X' +
H^{\dot{c}}_{\tilde{n}}=2^\omega$. 
On the other hand, $p^\star \geq p_{\xi_j}$, $j \leq K(2^{-\ell})$, so
$p^\star \forces_\lambda \dot{z} \not \in
X'+H^{\dot{c}}_{\tilde{n}}$. Contradiction.

To finish the proof we show that $\V^{{\mathcal P}_\lambda} \thinks 
\add(\M)=\kappa$.  
First note that ${\mathbf {MA}}_{<\kappa}$ implies that 
$\add(\M) \geq \kappa$ in $\V^{{\mathcal P}_\lambda}$.
The other inequality is a consequence of the general theory.
Recall that (see \cite{BJbook})
\begin{enumerate}
\item $\add(\M)=\min\{\cov(\M), {\mathfrak b}\}$
\end{enumerate}
Suppose that $F \subset  \omega^\omega$ is an unbounded family
of size $\geq \kappa$.
\begin{enumerate}
\item[2.] 
if $\mathcal P$ is a forcing notion of cardinality $<\kappa$ then 
$F$ remains unbounded in $\V^{\mathcal P}$.
\item[3.] if $\{{\mathcal P}_\alpha,{\mathcal Q}_\alpha: \alpha<\lambda\}$ is a
finite support iteration such that
$\forces_\alpha |{\mathcal Q}_\alpha|<\kappa$ then $\V^{{\mathcal
P}_\lambda} \thinks \text{$F$ is unbounded.}$.
\end{enumerate}

>From the results quoted above  follows that 
$\add(\M) \leq {\mathfrak b}\leq \kappa$ in $\V^{{\mathcal P}_\lambda}$,
which ends the proof.~$\QED$

\section{Strong measure zero sets}
In this section we will discuss models with strong measure zero sets
of size $ 2^{\boldsymbol\aleph_0} $.

We start with the definition of forcing that will be used in our construction.
\begin{definition}
The infinitely equal forcing notion $\mathbf {EE}$ 
 is defined as follows:
$p \in \mathbf {EE}$ if the following conditions are satisfied:
\begin{enumerate}
 \item $p:\dom(p) \longrightarrow 2^{<\omega}$,
  \item $\dom(p) \subseteq \omega, \ |\omega \setminus
    \dom(p)|=\boldsymbol\aleph_0$,
 \item $p(n) \in 2^n$ for all $n \in \dom(p)$.

\end{enumerate}
For $p,q \in \mathbf {EE}$ and $n \in \omega$ we define:
\begin{enumerate}
\item $p \geq q \iff p \supseteq q$, and 
\item $p \geq_n q \iff p \geq q$ and the first $n$ elements of $\omega
  \setminus \dom(p)$ and $\omega
  \setminus \dom(q)$ are the same.
\end{enumerate}
  
\end{definition}

It is easy to see (see \cite{BJbook}) that $\mathbf {EE}$
is proper (satisfies axiom A), and strongly $\omega^\omega $
bounding, that is if $p \forces \tau \in \omega $ and $ n \in \omega $
then there is $q \geq_n p$ and a finite set $F \subseteq \omega $ such
that $q \forces \tau \in F$.

In \cite{GJS92} it is shown that a countable support iteration of
$\mathbf {EE}$ and rational perfect set forcing produces a model where
there is a strong measure zero set of size $ 2^{\boldsymbol\aleph_0}
$. In particular, one can construct (consistently) a strong measure
zero of size $ 2^{\boldsymbol\aleph_0} $ without Cohen reals. 
The remaining question is whether such a construction can be carried out
without unbounded reals. 

\begin{theorem}[\cite{GJS92}]\label{omombdnosmz}
Suppose that $\{{\mathcal P}_\alpha, \dot{{\mathcal Q}}_\alpha:
\alpha<\omega_2\}$ is  a countable support iteration of proper, strongly
$\omega^\omega$-bounding forcing notions.
Then 
$$\V^{{\mathcal P}_{\omega_2}} \models \mathcal  {SN} \subseteq [\reals]^{\leq
  \boldsymbol \aleph_1}.~\QED$$
\end{theorem}
The theorem above shows that using countable support
iteration we cannot build a model with a strong measure zero set of size $>
{\mathfrak d}$. 
Since countable support iteration seems to be the universal method for
constructing  models with $ 2^{\boldsymbol\aleph_0}
=\boldsymbol\aleph_2$ the above result seems to indicate that a strong
measure zero set of size $> {\mathfrak d}$ cannot be constructed at all.
Strangely it is not the case. 

\begin{theorem}
  It is consistent that $\unif(\SN)=2^{\boldsymbol\aleph_0}>{\mathfrak
    d}=\boldsymbol\aleph_1 $ and there are strong measure zero sets of
  size $ 2^{\boldsymbol\aleph_0} $.
\end{theorem}
\Proof
Suppose that $\V \thinks \CH$ and $\kappa
=\kappa^{\boldsymbol\aleph_0}> \boldsymbol\aleph_1 $.
Let $ {\mathcal P} $ be a countable support product of $\kappa$ copies
of $\mathbf {EE}$.
The following facts are well-known (see \cite{GoShMan93})
\begin{enumerate}
\item $ {\mathcal P} $ is proper,
\item $ {\mathcal P} $ satisfies $ \boldsymbol\aleph_2 $-cc,
\item ${\mathcal P} $ is $\omega^\omega $-bounding,
\item for $f \in \V[G]\cap \omega^\omega $ there exists a countable
  set $A \subseteq \kappa $, $A \in \V$ such that $f \in \V[G \rest
  A]$.
\end{enumerate}
It follows from (3) that $\V^{{\mathcal P}} \thinks {\mathfrak
  d}=\boldsymbol\aleph_1 $. 
Moreover, (1) and (2) imply that $ 2^{\boldsymbol\aleph_0} =\kappa $
in $\V^{{\mathcal P}}$.

For a set $X \subseteq 2^\omega \cap \V^{\mathcal P} $ let $\supp(X)
\subseteq \kappa $ be a set such that $X \in \V[G \rest \supp(X)]$.
Note that $\supp(X)$ is not determined uniquely, but we can always
choose it so that $|\supp(X)|=|X|+\boldsymbol\aleph_0 $.

\begin{lemma}
  Suppose that $X \subseteq 2^\omega \cap \V^{\mathcal P}$ and $\supp(X) \neq \kappa $. Then $\V^{\mathcal P} \thinks
  X \in \SN$
\end{lemma}
Note that this lemma finishes the proof. Clearly the assumptions of
the lemma are met for all sets of size $<\kappa $ and also for many sets of
size $\kappa $.

\Proof
We will use the following characterization (see \cite{BJbook}):
\begin{lemma}\label{simplecharsmz}
The following conditions are equivalent.
\begin{enumerate}
\item $X \subseteq 2^\omega$ has strong measure zero.
\item For every $f \in \omega^\omega $ there exists $g \in
(2^{<\omega})^\omega$ such that $g(n) \in 2^{f(n)}$ for all $n$
and
$$\forall x \in X \ \exists n \ x \rest f(n) = g(n).~\QED$$
\end{enumerate}
\end{lemma}

Suppose that $X \subseteq \V^{\mathcal P} \cap 2^\omega $ is given and
$\supp(X)\neq \kappa $. Let $ \alpha^\star \in \kappa  \setminus
\supp(X)$. We will check condition (2) of the previous lemma. 

Fix $ f \in \V^{\mathcal P} \cap \omega^\omega $. Since ${\mathcal P}
$ is $\omega^\omega $-bounding we can assume that $f \in \V$.
Consider a condition $p \in {\mathcal P} $. Fix $\{k_n: n \in
\omega\}$ such that $k_n \geq f(n)$ and $k_n \not \in
\dom(p(\alpha^\star))$ for $n \in \omega $.
Let $p_f \geq p$ be any condition such that $\omega \setminus \{k_n :
n \in \omega\} \subseteq \dom(p_f(\alpha^\star))$.
We will check that
$$p_f \forces_{{\mathcal P}} \forall x \in X \ \exists n \ x \rest
f(n)=\dot{G}(\alpha^\star)(k_n) \rest f(n),$$
where $\dot{G}$ is the canonical name for the generic object. 
Take $x \in X$ and $r \geq p_f$.
Find $n $ such that $k_n \not \in \dom\lft1(r(\alpha^\star)\rgt1)$.
Let $r' \geq r$ and $s$ be such that 
\begin{enumerate}
\item $\supp(r') \subseteq \supp(X)$
\item $r' \geq r \rest \supp(X)$,
\item $r' \forces_{{\mathcal P}} x \rest k_n=s$.
\end{enumerate}
Let 
$$r''(\beta)=\left\{
  \begin{array}{ll}
r'(\beta)& \text{if } \beta\neq \alpha^\star\\
r'(\alpha^\star)\cup \{(k_n,s)\}& \text{if }\beta=\alpha^\star
  \end{array}\right. .$$
It is easy to see that $r''\forces x \rest f(n) =
\dot{G}(\alpha^\star)(k_n) \rest f(n)$.
Since  $f$ and $x$  were arbitrary we are done.~$\QED$

\end{document}